\documentclass{article}
\usepackage{amsmath}

\newtheorem{theorem}{Theorem}

\newtheorem{proposition}[theorem]{Proposition}

\input{tcilatex}
\begin{document}

\title{Affine Connection Induced from The Horizontal lift $^{H}\nabla $ on a
Cross-section}
\author{Melek ARAS\thanks{%
Department of mathematics, Faculty of Art and Sciences, Giresun
Universty,28049 Giresun, Turkey,
E-mail:melekaras25@hotmail.com:melek.aras@giresun.edu.tr}}
\maketitle

\begin{abstract}
The main purpose of present paper is to study the affine connection induced
\ from the horizontal lift $\overline{\bigtriangledown }$ on the
cross-section $\beta _{\vartheta }\left( M_{n}\right) $ determined \ by a
vector field $\vartheta $ in $M_{n}$ with respect to the adapte frame of $%
\beta _{\vartheta }\left( M_{n}\right) $.


\textit{2010 AMS Classification:53C05, 53B05, 53C07}
\end{abstract}

\begin{description}
\item[1. Introduction] 
\end{description}

Let $M_{n}$ be an \textit{n-}dimensional differentiable manifold of class $%
C^{\infty }$ an $T_{p}\left( M_{n}\right) $ the tangent space at a point $P$
of $M_{n}$, that is, the set of all tangent vectors of $M_{n}$ at $P$. Then
the set

\begin{equation*}
T\left( M_{n}\right) =\underset{P\epsilon M_{n}}{\cup }T_{P}\left(
M_{n}\right) ,
\end{equation*}%
is by definition, tangent bundle over the manifold $M_{n}$ \cite{1}.

Let $M_{n}$ be a Riemannian manifold with metric $g$ whose components in a
coordinate neighborhood $U$ are $g_{ji}$, and denote by $\Gamma _{ij}^{k}$
the Christoffel symbols formed with $g_{ji}.$If $U$ being a neighborhood of $%
M_{n}$, then the horizontal lift $^{H}g$ of $g$ has components

\begin{equation*}
^{H}g=\left( 
\begin{array}{cc}
\Gamma _{i}^{m}g_{mj}+\Gamma _{j}^{m}g_{im} & g_{ij} \\ 
g_{ij} & 0%
\end{array}%
\right)
\end{equation*}

with respect to the induced coordinates $\left( x^{h},y^{h}\right) $ in $\pi
^{-1}\left( U\right) \subset T\left( M_{n}\right) $, where $\Gamma
_{i}^{m}=y^{j}\Gamma _{ji}^{m}$, $\Gamma _{ji}^{m}$ being the components of
the affine connection in $M_{n}.$

Now we shall define the horizontal lift $\overline{\nabla }$ of the affine
connection $\nabla $ in $M_{n}$ to $T\left( M_{n}\right) $ by the conditions

\begin{equation}
\begin{array}{l}
\overline{\nabla }_{^{V}X}\text{ }^{V}Y=0\text{,\ \ \ \ \ \ \ \ \ \ \ \ \ \
\ \ \ }\overline{\nabla }_{^{V}X}\text{ }^{H}Y=0 \\ 
\multicolumn{1}{c}{\overline{\nabla }_{^{H}X}\text{ }^{V}Y=^{V}\left( \nabla
_{X}Y\right) \text{, \ \ \ \ \ \ }\overline{\nabla }_{^{H}X}\text{ }%
^{H}Y=^{H}\left( \nabla _{X}Y\right) ,\text{\ \ \ \ \ \ \ \ \ \ \ \ \ \ \ \
\ \ \ }} \\ 
\multicolumn{1}{c}{}%
\end{array}
\tag{1}
\end{equation}

for $X,Y\epsilon \Im _{0}^{1}\left( M_{n}\right) $. From $\left( 1\right) $,
the horizontal lift $\overline{\nabla }$ of $\nabla $ has components $%
\overline{\Gamma }_{JI}^{K}$ such that

\begin{equation}
\begin{array}{c}
\overline{\Gamma }_{ji}^{k}=\Gamma _{ji}^{k},\text{\ \ \ \ \ \ }\overline{%
\Gamma }\text{\ }_{j\overline{i}}^{k}=\overline{\Gamma }_{\overline{j}i}^{k}=%
\overline{\Gamma }_{\overline{j}\overline{i}}^{k}=\overline{\Gamma }_{%
\overline{j}\overline{i}}^{\overline{k}}=0, \\ 
\overline{\Gamma }_{ji}^{\overline{k}}=y^{s}\partial _{s}\Gamma
_{ji}^{k}-y^{s}R_{sji}^{k},\text{ \ \ \ \ }\overline{\Gamma }_{\overline{j}%
i}^{\overline{k}}=\overline{\Gamma }_{j\overline{i}}^{\overline{k}}=\Gamma
_{ji}^{k}%
\end{array}
\tag{2}
\end{equation}%
\ \ with respect to the induced coordinates in $T\left( M_{n}\right) ,$where 
$\Gamma _{ji}^{k}$ are the components of $\nabla $ in $M_{n}$ \cite{6}.

Let a vector-field in a manifold $M_{n}$, then the vector field defines a
cross-section in the tangent bundle $T\left( M_{n}\right) $. Tensor fields
and Connections on a Cross-Section in the Tangent bundel was studied by Houh
and Ishihara $\cite{1}$, Tani $\cite{3},$Yano $\cite{4}$. Affine connections
induced from $\nabla ^{C}$ on the cross-section $\beta _{\vartheta }\left(
M_{n}\right) $ was studied by Yano and Ishihara \cite{6}.

We suppose that there is give a vector field $\vartheta $ in an \textit{n}%
-dimensional manifold $M_{n}$. Then the correspondence $p\rightarrow
\vartheta _{p}$, $\vartheta _{p}$ being the value of $\vartheta $ at $p\in
M_{n}$, determines a mapping $\beta _{\vartheta }:M_{n}\rightarrow T\left(
M_{n}\right) $ and the \textit{n}-dimensional submanifold $\beta _{\vartheta
}\left( M_{n}\right) $ of $T\left( M_{n}\right) $ is called the
cross-section determined by $\vartheta $. If the vector field $\vartheta $
has local components $\vartheta ^{k}\left( x\right) $ in $M_{n}$. Then the
cross-section $\beta _{\vartheta }\left( M_{n}\right) $ is locally expressed
by

\begin{equation}
x^{h}=x^{h}\text{ , \ \ \ \ \ \ \ \ \ \ \ }y^{h}=\vartheta ^{h}\left(
x\right)  \tag{3}
\end{equation}

with respect to the induced coordinates $\left( x^{A}\right) =\left(
x^{h},y^{h}\right) $ in $T\left( M_{n}\right) $. Differentiating $\left(
3\right) $, we see that n tangent vectors $B_{\left( j\right) }$ to $\beta
_{\vartheta }\left( M_{n}\right) $ have components

\begin{equation*}
B_{j}^{A}=\frac{\partial x^{A}}{\partial x^{j}}
\end{equation*}

i.e.,

\begin{equation}
B_{\left( j\right) }:\left( B_{j}^{A}\right) =\left( 
\begin{array}{c}
\delta _{j}^{h} \\ 
\partial _{j}\vartheta ^{h}%
\end{array}%
\right)  \tag{4}
\end{equation}%
with respect to the induced coordinates $T\left( M_{n}\right) $.

On the other hand, since a fibre is locally expressed by $x^{h}=const.$, $%
y^{h}=y^{h}$, $y^{h}$ being considered as parameters.

\begin{equation}
C_{\left( j\right) }:\left( C_{\left( j\right) }^{A}\right) =\left( 
\begin{array}{c}
0 \\ 
\delta _{j}^{h}%
\end{array}%
\right)  \tag{5}
\end{equation}%
are tangent to the fibre.

We now consider in $\pi ^{-1}\left( U\right) $ , $U$ being coordinate
neighborhood of $M_{n}$, $2n$ local vector fields $B_{\left( j\right) }$ and 
$C_{\left( j\right) }$ along $\beta _{\vartheta }\left( M_{n}\right) $,
represented respectively by

\begin{equation*}
B_{\left( j\right) }=B\frac{\partial }{\partial x^{j}}\text{, \ \ \ \ \ \ \ }%
C_{\left( j\right) }=C\frac{\partial }{\partial x^{j}}\text{\ .\ }
\end{equation*}

They form a local family of frames $\left\{ B_{\left( j\right) },C_{\left(
j\right) }\right\} $ along $\beta _{\vartheta }\left( M_{n}\right) $, which
is called the adapted frame of $\beta _{\vartheta }\left( M_{n}\right) $ in $%
\pi ^{-1}\left( U\right) $ $\cite{6}$

\begin{description}
\item[2.Affine Connection Induced from $\overline{\protect\nabla }$ on a
Cross-Section] 
\end{description}

We suppose that $M_{n}$ is a manifold with affine connection $\nabla $. Thus
the tangent bundle $T\left( M_{n}\right) $ of $M_{n}$ is a manifold with
affine connection $\overline{\nabla }$ which is the horizontal lift of $%
\nabla $. We now study the affine connection induced from $\overline{\nabla }
$ on the cross-section $\beta _{\vartheta }\left( M_{n}\right) $ determined
by a vector field $\vartheta $ in $M_{n}$ with respect to the adapted frame
of $\beta _{\vartheta }\left( M_{n}\right) .$

The lineer connection $^{^{\prime }}\nabla $ on the cross-section $\beta
_{\vartheta }\left( M_{n}\right) $ induced from $\overline{\nabla }$ is
defined by connection components $^{^{\prime }}\Gamma _{ji}^{h}$ given by 
\cite{6}

\begin{equation}
^{^{\prime }}\Gamma _{ji}^{h}=\left( \partial _{j}B_{i}^{A}+\overline{\Gamma 
}_{MN}^{A}B_{j}^{M}B_{i}^{N}\right) B_{A}^{h},  \tag{6}
\end{equation}%
where $\overline{\Gamma }_{MN}^{A}$ are the connection components of $%
\overline{\nabla }$ with respect to the induced coordinates in $T\left(
M_{n}\right) $ and $B_{A}^{h}$ are defined by

\begin{equation*}
\left( B_{A}^{h},C_{A}^{h}\right) =\left( B_{j}^{A},C_{j}^{A}\right) ^{-1}
\end{equation*}%
and hence 
\begin{equation}
\left( B_{B}^{h}\right) =\left( \delta _{j}^{h},0\right) ,\text{ \ \ \ \ \ \
\ \ }\left( C_{B}^{h}\right) =\left( -\partial _{j}\vartheta ^{h},\delta
_{j}^{h}\right) .  \tag{7}
\end{equation}

Substituting $\left( 2\right) $ for $\overline{\Gamma }_{MN}^{A}$, $\left(
4\right) $, $\left( 5\right) $ and $\left( 7\right) $ in $\left( 6\right) $,
we find

\begin{equation}
^{^{\prime }}\Gamma _{ji}^{h}=\Gamma _{ji}^{h}  \tag{8}
\end{equation}%
where $\Gamma _{ji}^{h}$ are components of $\nabla $ in $M_{n}$.

From $\left( 6\right) $ we see that

\begin{equation}
\partial _{j}B_{i}^{A}+\overline{\Gamma }_{MN}^{A}B_{j}^{M}B_{i}^{N}-\Gamma
_{ji}^{h}B_{h}^{A}=H_{ji}^{\overline{k}}C_{\left( k\right) }^{A},  \tag{9}
\end{equation}%
i.e., that the left hand side is a linear combinations of $C_{\left(
h\right) }^{A},$ where the coefficients $H_{ji}^{\overline{h}}$ \ will be
found in the sequel. To find the coefficients $H_{ji}^{\overline{h}},$ we
put $A=h$ in $\left( 9\right) $ and hence obtain

\begin{equation}
H_{ji}^{\overline{h}}=\partial _{j}\partial _{i}\vartheta ^{h}+\vartheta
^{t}\partial _{t}\Gamma _{ji}^{h}+\vartheta ^{t}R_{tji}^{h}+\Gamma
_{mi}^{h}\partial _{j}\vartheta ^{m}+\Gamma _{jn}^{h}\partial _{i}\vartheta
^{n}-\Gamma _{ji}^{t}\partial _{t}\vartheta ^{h}  \tag{10}
\end{equation}%
which are components $\mathit{\tciLaplace }_{\vartheta }\Gamma _{ji}^{k}$ of
the Lie derivative of the affine connection $\nabla $ with respect to $%
\vartheta $\cite{5}. Thuse, representing the left-hand side of $\left(
9\right) $ by $^{^{\prime }}\nabla _{j}B_{i}^{A}$, we have from $\left(
10\right) $

\begin{equation}
^{^{\prime }}\nabla _{j}B_{i}^{A}=\left( \mathit{\tciLaplace }_{\vartheta
}\Gamma _{ji}^{h}+\vartheta ^{t}R_{tji}^{h}\right) C_{h}^{A}.  \tag{11}
\end{equation}%
Thus we have

\begin{proposition}
If $\vartheta ^{t}R_{tji}^{h}=0$, then $^{^{\prime }}\nabla
_{j}B_{i}^{A}=\left( \mathit{\tciLaplace }_{\vartheta }\Gamma
_{ji}^{h}+\vartheta ^{t}R_{tji}^{h}\right) C_{h}^{A}$ is the equation of
Gauss for the cross-section $\beta _{\vartheta }\left( M_{n}\right) $
determined by a vector field $\vartheta $ in $M_{n}$ to $T\left(
M_{n}\right) .$
\end{proposition}

\begin{proposition}
In order that the cross-section in $T\left( M_{n}\right) $ determined by a
vector field $\vartheta $ in $M_{n}$ with affine connection $\nabla $ be
totally geodesic with respect to $\overline{\nabla }$ it is necessary and \
sufficient \ that respectively $\vartheta $ is an infinitesimal \ affine
transformation in $M_{n}$, i.e.,that $\tciLaplace _{\vartheta }\nabla =0$ 
\cite{4} \ and $\vartheta ^{t}R_{tji}^{h}=0,$where $R_{tji}^{h}$ is
components of the curvature tensor $R$ of $\nabla $.
\end{proposition}

By means of $\left( 9\right) $, the equation $\left( 11\right) $ reduces to

\begin{equation}
\overline{\nabla }_{B_{\left( j\right) }}B_{\left( i\right) }=\Gamma
_{ji}^{h}B_{\left( h\right) }+H_{ji}^{\overline{h}}C_{\left( h\right) }. 
\tag{12}
\end{equation}

We now have

\begin{equation}
\overline{R}\left( B_{\left( k\right) },B_{\left( j\right) }\right)
B_{\left( i\right) }=\overline{\nabla }_{B_{\left( k\right) }}\overline{%
\nabla }_{B_{\left( j\right) }}B_{\left( i\right) }-\overline{\nabla }%
_{B_{\left( j\right) }}\overline{\nabla }_{B_{\left( k\right) }}B_{\left(
i\right) },  \tag{13}
\end{equation}

$\overline{R}$ being the curvature tensor of $\overline{\nabla }$ because $%
\left[ B_{\left( j\right) },B_{\left( i\right) }\right] =0$. Thus, denoting
by $R_{kji}^{h}B_{\left( h\right) }$ the components of the curvature tensor $%
R$ \ of $\nabla $, we have from $\left( 13\right) $

\begin{equation}
\begin{array}{c}
\overline{R}\left( B_{\left( k\right) },B_{\left( j\right) }\right)
B_{\left( i\right) }=R_{kji}^{h}B_{\left( h\right) }+\left\{ \nabla
_{k}\left( \mathit{\tciLaplace }_{\vartheta }\Gamma _{ji}^{h}\right) -\nabla
_{j}\left( \mathit{\tciLaplace }_{\vartheta }\Gamma _{ki}^{h}\right)
\right\} C_{\left( h\right) } \\ 
+\nabla _{k}\left( \mathit{\vartheta }^{t}R_{tji}^{h}\right) -\nabla
_{j}\left( \mathit{\vartheta }^{t}R_{tki}^{h}\right) \\ 
\end{array}
\tag{14}
\end{equation}

which reduces to

\begin{equation}
\overline{R}\left( B_{\left( k\right) },B_{\left( j\right) }\right)
B_{\left( i\right) }=R_{kji}^{h}B_{\left( h\right) }+\left( \mathit{%
\tciLaplace }_{\vartheta }R_{kji}^{h}\right) C_{\left( h\right) }+\nabla
_{k}\left( \mathit{\vartheta }^{t}R_{tji}^{h}\right) -\nabla _{j}\left( 
\mathit{\vartheta }^{t}R_{tki}^{h}\right)  \tag{15}
\end{equation}

where the well know formula \cite{5}

\begin{equation*}
\nabla _{k}\left( \mathit{\tciLaplace }_{\vartheta }\Gamma _{ji}^{h}\right)
-\nabla _{j}\left( \mathit{\tciLaplace }_{\vartheta }\Gamma _{ki}^{h}\right)
=\mathit{\tciLaplace }_{\vartheta }R_{kji}^{h}.
\end{equation*}

from $\left( 15\right) $, we have

\begin{proposition}
In order that $\overline{R}\left( \widetilde{X},\widetilde{Y}\right) 
\widetilde{Z}$ evaluated for vector fields $\widetilde{X},\widetilde{Y}$ and 
$\widetilde{Z}$ tangent to the cross-section determined by a vector field $%
\vartheta $ in $M_{n},R$ being curvature tensor of an affine connection $%
\nabla $, be always tangent to the cross-section, it is necessary and
sufficient that respectively the Lie derivative $\mathit{\tciLaplace }%
_{\vartheta }R$ of $R$ with respect to $\vartheta $ in $M_{n}$ vanishes,
i.e., $\mathit{\tciLaplace }_{\vartheta }R=0\cite{4}$ and $\nabla _{k}\left( 
\mathit{\vartheta }^{t}R_{tji}^{h}\right) -\nabla _{j}\left( \mathit{%
\vartheta }^{t}R_{tki}^{h}\right) =0.$
\end{proposition}

\end{document}